\documentclass[11pt, 2side]{amsart}    

\usepackage{graphicx}
\usepackage{amssymb,amsmath}

\numberwithin{equation}{section}
\numberwithin{figure}{section}%
\newtheorem{theo}{Theorem}

\newtheorem{coro}{Corollary}

\newtheorem{lema}{Lemma}

\newcommand{\ba}{\begin{eqnarray} }
\newcommand{\ea}{\end{eqnarray} }
\newcommand{\be}{\begin{equation} }
\newcommand{\ee}{\end{equation} }
\newenvironment{dwd}{\par\noindent{\bf Proof.}}{\par\rightline{$\blacksquare$}}

\newcommand{\ra}{\rightarrow}

\newcommand{\E}{{\mathbf E}}
\newcommand{\PM}{{\mathbf P}}

\newcommand{\ls}{\leq}
\newcommand{\gs}{\geq}
\begin{document}

\title[Truncated variation] {Moments approach to the concentration properties of truncated variation}

\author{Witold Bednorz         \and
        Rafa\l{} \L{}ochowski
}

\address{Institute of Mathematics, Warsaw University, Banacha 2, 02-097 	Warszawa, Poland}
\email{wbednorz@mimuw.edu.pl}           
\address{Department of Mathematics and Mathematical Economics, Warsaw School of Economics, Madali\'{n}skiego 6/8, 02-513 Warszawa, Poland
}
\email{rlocho@sgh.waw.pl}

\subjclass[2010]{60G17, 60G15, 60G22, 60G48, 60G51}
\keywords{sample boundedness, Gaussian processes, diffusion, L\'{e}vy processes, truncated variation}

\thanks{ The first author's research is supported by the NCN grant DEC-2012/05/B/ST1/00412.}

\begin{abstract}
In this paper we show an alternative approach to the concentration of truncated variation
for stochastic processes on a real line. Our method is based on the moments control and can be used to
generalize the results to the case of processes with heavy tails. 
\end{abstract}

\maketitle

\section{Introduction}

Let $X=X(t)$, $t\in [0,1]$ be a real valued
stochastic process with c\`{a}dl\`{a}g trajectories. 
In this paper we are interested in concentration properties of
the truncated variation of $X$. It is well-known that usually the total variation
for semimartingales is infinite, whereas it is always possible to
define for a given $c>0$  
\begin{equation}
TV^{c}\left(X\right) :=\sup_{n}\sup_{0=
t_{0}<t_{1}<...<t_{n}= 1}\sum_{i=1}^{n}\max \left\{ \left| X(
t_{i}) -X(t_{i-1}) \right| -c,0\right\}.  \label{tv:def}
\end{equation}
It was proved in \cite{LM} that for any continuous semimartingale $X$ on $[0,1]$ we have  
\begin{equation} \label{LLN}
c \cdot TV^{c}\left(X\right) \rightarrow_{c\downarrow 0} 
\langle X\rangle_1 -\langle X \rangle_{0} \text{ a.s.,}
\end{equation}
where $\left\langle .\right\rangle $ denotes the quadratic variation
of $X$. In this way truncated variation may be useful in the so called pathwise approach to stochastic integration. 
\smallskip

\noindent
In the paper \cite{Be-Lo} we have proved concentration inequalities for various processes whose 
increments decay exponentially fast. In this paper we introduce a new approach based on the method of moments. In this way it is possible to establish bounds on moments of the truncated variation under 
much weaker assumptions. In particular it should be possible to study some processes with heavy tailed
increments.

\section{The moments control}\label{sect2}

We assume that for a given $k\gs 1$ and all $s,t\in [0,1]$
\be\label{nier:0}
\|X(t)-X(s)\|_k\ls C_1k^{p}|s-t|^q,
\ee
where $0<q<1$, $p>0$. The second condition we require is that if $d\gs C_2k^p|s-t|^q$ for some
$s,t\in [0,1]$ then also 
\be\label{nier:1}
\|(|X(t)-X(s)|-d)_{+}\|_{k}\ls C_2 k^{p}|s-t|^q f(\frac{d}{C_2k^p|s-t|^q}),
\ee
where the function $f$ is positive, decreasing and satisfies the following growth condition: for some integer $r\gs 2^{1/q}$ and any given constant $C_3$ the following must holds
\be\label{nier:2}
\sum^{\infty}_{l=0}r^{l(1-q)}f(C_3(r^q/2)^{l/p})<\infty.
\ee
It means that the function $f$ must decay slightly faster than the function $1/x^{\alpha}$,
where $\alpha>\frac{p}{q}-p$.
The natural setting in which these condition works is when (\ref{nier:0})
is satisfied for all $k\gs 1$ and $p\ls 1$. Then obviously
$$
\PM(|X(t)-X(s)|>xeC_1|s-t|^q)\ls \exp(-x^{\frac{1}{p}})\;\;\mbox{for}\;x\gs 1
$$
and thus for all $x>eC_1|s-t|^q$
$$
\PM(|X(t)-X(s)|>x)\ls \exp(-\left[\frac{x}{(eC_1|s-t|^q)}\right]^{1/p}).
$$
Thus if $d>eC_1|s-t|^q$ then
\begin{align*}
& \E|(|X(t)-X(s)|-d)_{+}|^k=\int^{\infty}_0 kx^{k-1}\PM(|X(t)-X(s)|>x+d)dx\\
& =\int^{\infty}_0 kx^{k-1}\exp(-\left[\frac{d+x}{eC_1|s-t|^q}\right]^{1/p}).
\end{align*}
Using that $p\ls 1$ and $x>0$ we obtain
$$
\left[\frac{d+x}{eC_1|s-t|^q}\right]^{1/p}\gs 
\left[\frac{d}{eC_1|s-t|^q}\right]^{1/p}+\left[\frac{x}{eC_1|s-t|^q}\right]^{1/p}
$$
and hence
\begin{align*}
& \int^{\infty}_0 kx^{k-1}\exp(-\left[\frac{d+x}{eC_1|s-t|^q}\right]^{1/p}) \\
& \ls\exp(-\left[\frac{d}{eC_1|s-t|^q}\right]^{\frac{1}{p}})(eC_1|s-t|^q)^k\int^{\infty}_{0} kp x^{kp-1}\exp(-x)dx\\
&= \exp(-\left[\frac{d}{eC_1|s-t|^q}\right]^{\frac{1}{p}})(eC_1|s-t|^q)^k \Gamma(kp).
\end{align*}
Since $[\Gamma(kp)]^{\frac{1}{k}}\ls (kp)^p$ we finally get
$$
\|(X(t)-X(s))-d)_{+}\|_{k}\ls (kp)^p(eC_1|s-t|^q)\exp(-\frac{1}{k}\left[\frac{d}{eC_1|s-t|^q}\right]^{\frac{1}{p}}).
$$
Therefore (\ref{nier:1}) holds with $C_2=eC_1p^p$ and  
We use (\ref{nier:0}) and (\ref{nier:1}) to control moments of the truncated variation.
\begin{theo}\label{thm1}
Suppose that (\ref{nier:0}) and (\ref{nier:1}) are satisfied.
Then there exists a universal constant $K$ which may depend on $p,q$ such that  
$$
\|TV^c(X,T)\|_k\ls Kc^{1-\frac{1}{q}}k^{\frac{p}{q}}.
$$
\end{theo}
\begin{dwd}
Consider given partition $\Pi_d=\{t_0,t_1,\dots,t_d\}$, where $0\ls t_0<t_1<\dots<t_d\ls 1$. 
We have to provide a universal upper bound for
$$
\sum^d_{i=1}\left(|X(t_i)-X(t_{i-1})|-c\right)_{+}.
$$
In order to find such an estimate following paper \cite{Be-Lo} we introduce
the approximation sequence $(T_n)^{\infty}_{n=0}$, $T_n\subset T$.
Fix integer $r>1$ and let $T_n=\{kr^{-n}:\;k=0,1,\ldots,r^{n}\}$.
Obviously $T_{n}\subset T_{n+1}$ and $|T_n|=r^n+1$. We define neighbourhood
of a given point $t\in T_{n+1}\supset T_n$, namely
$$
I_{n+1}(t)=\{s\in T_{n+1}: |s-t|< 2r^{-n}\}.
$$
Clearly $|I_{n+1}(t)|\ls 2r-1$. 
For a given set $T_n$  and $t \in T$, by $\pi_n(t)$ we denote the unique point $s\in T_n$ such that $s\ls t$ and $|t-s|< r^{-n}$.  This way we define the function $\pi_n:T\ra T_n$.
We have $|t-\pi_n(t)| < r^{-n}$ for all $t\in T$ 
and $\pi_n(s)\ls \pi_n(t)$ if $s\ls t$. Note that $s,t\in T_n$, $s\neq t$,
$|s-t|\gs r^{-n}$. We use the above construction to approximate intervals
$[t_{i-1},t_{i}]$ for any consecutive $t_{i-1},t_{i}\in \Pi_d$ where $i\in \{1,2,\ldots,d\}$.
We denote by $t^n_i=\pi_n(t_i)$, it is crucial to observe that for a fixed $t^n_i$
there cannot be to many candidates for $t^{n+1}_i$. Clearly 
$$
|t^n_i-t^{n+1}_i|\ls |t_i-t^n_i|+|t_i-t^{n+1}_i|< r^{-n}+r^{-n-1}<2r^{-n}.
$$
Consequently $t^{n+1}_i\in I_{n+1}(t^n_i)$. Moreover 
$t_0^{n+1}\ls t_1^{n+1}\ls \ldots \ls t_d^{n+1}$. Obviously $(t_i^n)^{\infty}_{n=0}$
is a path that approximates point $t_i$ in the sense that $\lim_{n\ra \infty}|t_i-t^n_{i}|=0$.
Since we have to consider all intervals $[t_{i-1},t_i]$, $i=1,2,\ldots,d$ we first classify
their lengths. To this aim we define for $m=0,1,\ldots$
$$
J_m=\{i\in\{1,\ldots,d\}:\;r^{-m-1}<|t_{i-1}-t_i|\ls r^{-m}\}.
$$
The approximation paths for the interval $[t_{i-1},t_i]$ 
consists of $(t^n_i)^{\infty}_{n=m+1}$ and $(t^n_{i-1})^{\infty}_{n=m+1}$.
We call any pair $(t^n_{i},t^{n+1}_i)$, $n\gs m+1$ a step of the approxiamtion.
Observe that 
$$
|t^{m+1}_i-t^{m+1}_{i-1}|\ls |t_i-t^{m+1}_i|+|t_i-t_{i-1}|+|t_{i-1}-t^{m+1}_{i-1}|<
2r^{-m-1}+r^{-m}\ls 2r^{-m}.
$$
Therefore $t^{m+1}_{i-1}\in I_{m+1}(t^{m+1}_i)$ and $t^{m+1}_i\in I_{m+1}(t^m_{i-1})$.
We have to prove that for fixed $u\in T_n$, $v\in I_{n+1}(u)\subset T_{n+1}$
there are at most two different intervals $[t_{i-1},t_i]$, $i=1,2,\ldots,d$ such that  
$t^n_i=u$, $t^{n+1}_i=v$ or $t^n_{i-1}=u$, $t^{n+1}_{i-1}=v$.
\begin{lema}\label{nowy1}
Consider $u\in T_{n}$, and $v\in I_{n+1}(u)\subset T_{n+1}$. The step $(u,v)$ may occur in 
the approximation of $\Pi_d$ in two ways: either there exists no more than one
$i\in \{1,\ldots,d\}$ such that $i\in J_m$, $m+1\ls n$ and $t^n_{i-1}=u\in T_n$, $t^{n+1}_{i-1}=v$ or
or there exists no more than one $i\in \{1,2,\ldots,d\}$ such that
$i\in J_{m'}$, $m'+1\ls n$ and $t^n_i=u,$ $t^{n+1}_{i}=v$. 
\end{lema}
\begin{dwd}
Recall that $r\gs 2$. It suffices to prove that for a given $i\in J_m$, $n\gs m+1$ points $t^{n+1}_i$ and $t^{n+1}_{i-1}$ are different. Indeed since $t^{n+1}_0\ls t^{n+1}_1\ls \dots\ls t^{n+1}_d$ the property implies that there can be at most one $i\in\{0,1,\dots,d\}$
such that $t^{n+1}_i=v$. To prove the assertion we use $|t_i-t_{i-1}|>r^{-m-1}$ which implies that for  $l \gs m+1$
\begin{align*}
&|t^{n+1}_i-t^{n+1}_{i-1}|\geq r^{-m-1}-|t^{n+1}_{i-1}-t_{i-1}|-|t^{n+1}_i-t_i| \\
&> r^{-m-1}-2r^{-n-1}\gs 0
\end{align*}
since $r\gs 2$.
\end{dwd}
Let $m_k$ be defined as $k^p r^{-(m_k+1)q}< c/M_0\ls k^p r^{-m_k q}$, where $M_0=2C_2er^{2q}$.
We are ready to state the right upper bound on the truncated variation.
\begin{lema}
The following estimate holds
\begin{align*}
&\sum^d_{i=1}\left(|X(t_i)-X(t_{i-1})|-c\right)_{+}\ls
2\sum^{m_k}_{n=0} \sum_{u\in T_{n+1}}\sum_{v\in I_{n+1}(u)}|X(u)-X(v)|\\
&+2\sum^{\infty}_{n=m_k+1} \sum_{u\in T_{n+1}}\sum_{v\in I_{n+1}(u)}(|X(u)-X(v)|-2^{m_k-n-1}c)_{+}
\end{align*}
\end{lema}
\begin{dwd}
First consider given interval $[t_{i-1},t_i]$ and $i\in J_m$. If $m\ls m_k$
\begin{align*}
&\left(|X(t_i)-X(t_{i-1})|-c\right)_{+}\ls |X(t^{m+1}_i)-X(t^{m+1}_{i-1})|+\\
&+\sum_{s\in\{i-1,i\}}
\sum^{m_k}_{n=m+1} |X(t^{n}_s)-X(t^{n+1}_s)|\\
&+\sum^{\infty}_{n=m_k+1}
(|X(t^n_s)-X(t^{n+1}_s)|-2^{m_k-n-1}c)_{+}.
\end{align*}
On the other hand if $i\in J_m$, $m>m_k$ then
\begin{align*}
&\left(|X(t_i)-X(t_{i-1})|-c\right)_{+}\\
&\ls (|X(t^{m+1}_i)-X(t^{m+1}_{i-1})|-2^{m_k-m-2}c)_{+}\\
&+\sum_{s\in\{i-1,i\}}\sum^{\infty}_{n=m+1}
(|X(t^n_s)-X(t^{n+1}_s)|-2^{m_k-n-1}c)_{+}.
\end{align*}
It suffices to apply Lemma \ref{nowy1} to finish the proof.
\end{dwd}
Consequently we can bound
\begin{align*}
&\|TV^c(X,T)\|_k\ls 2\sum^{m_k}_{n=0}\sum_{u\in T_{n+1}}\sum_{v\in I_{n+1}(u)}\|X(u)-X(v)\|_k\\
&+2\sum^{\infty}_{n=m_k+1} \sum_{u\in T_{n+1}}\sum_{v\in I_{n+1}(u)}
\|(|X(u)-X(v)|-2^{m_k-n-1}c)_{+}\|_k.
\end{align*}
By the assumption (\ref{nier:0}) we have
$$
\|X(u)-X(v)\|_k\ls C_1k^p|u-v|^q
$$
and hence
\begin{align*}
&2\sum^{m_k}_{n=0}\sum_{u\in T_{n+1}}\sum_{v\in I_{n+1}(u)}\|X(u)-X(v)\|_k
&\ls 
2C_1\sum^{m_k}_{n=0}\sum_{u\in T_{n+1}}\sum_{v\in I_{n+1}(u)}k^p|u-v|^q\ls \\
&\ls 2C_1\sum^{m_k}_{n=0}(r^{n+1}+1)(2r-1)k^p(2r^{-n})^q\\
&\ls 8C_1r^2k^p\sum^{m_k}_{n=0}r^{n(1-q)}.
\end{align*}
Using that $q<1$ and $r^{m_k(1-q)}\ls M_0^{\frac{1-q}{q}}k^{\frac{p}{q}-p} c^{-\frac{1-q}{q}}$ we establish the final bound for this part
\begin{align*}
& 2\sum^{m_k}_{n=0}\sum_{u\in T_{n+1}}\sum_{v\in I_{n+1}(u)}\|X(u)-X(v)\|_k \\
& \ls 8Cr^2\frac{k^p(r^{(m_k+1)(1-q)}-1)}{r^{1-q}-1} \ls D_1(r,q)k^{\frac{p}{q}}c^{1-\frac{1}{q}},
\end{align*}
where $D_1(r,q)=(8Cr^{3-q}M_0^{\frac{1}{q}-1})/(r^{1-q}-1)$.
Now we have to bound
$$
\|(|X(u)-X(v)|-2^{m_k-n-1}c)_{+}\|_k.
$$
For this part we need the second assumption (\ref{nier:1}). In order to use the inequality we need that
$d=2^{m_k-n-1}c\gs C_2k^p|u-v|^q$. Note that if $u\in T_{n+1}$, $v\in I_{n+1}(u)$, $n>m_k$
then $|u-v|\ls 2r^{-n}$ and hence
\begin{align*}
&C_2k^p|u-v|^q\ls 2^qC_2 k^p r^{-nq}\\
&\ls  2^qC_2r^{2q}k^p r^{-(n+1-m_k)q}r^{-(m_k+1)q}\ls 2^{-(n+1-m_k)q}c
\end{align*}
since $r\gs 2$ and $r^{-(m_k+1)q}<c/M_0$ and $M_0=2^qC_2r^{2q}$. Therefore by (\ref{nier:1}) 
$$
\|(|X(u)-X(v)|-d)_{+}\|_k\ls C_2 k^p|u-v|^q f(\frac{d}{C_2k^p|u-v|^q}).
$$
Note that for $u\in T_{n+1}$ and $v\in I_{n+1}(u)$ due to $d=2^{m_k-n-1}$
and $k^pc^{-1}\ls M_0^{-1}r^{(m_k+1)q}$ we get
\begin{align*}
& \frac{k^p|u-v|^q}{d}\ls 2^{n+1-m_k}(2r^{-n})^q k^pc^{-1}\ls\\ 
&\ls M_0^{-1}(2r^{-n})^q(2^{n+1-m_k}r^{q(m_k+1)})\ls 
4\cdot 2^q M_0^{-1}(r^q/2)^{-(n-m_k-1)}.
\end{align*}
Therefore
\begin{align}
&\|(|X(u)-X(v)|-2^{m_k-n-1}c)_{+}\|_k\\
\label{mela}&\ls C_2k^p r^{-qn}f(D_2(r,q)(r^q/2)^{-(n-m_k-1)/p}), 
\end{align}
where $D_2(r,q)=(4\cdot 2^q)^{-1} M_0$. It remains to sum up all the bounds (\ref{mela})
\begin{align*}
&2\sum^{\infty}_{n=m_k+1} \sum_{u\in T_{n+1}}\sum_{v\in I_{n+1}(u)}(|X(u)-X(v)|-2^{m_k-n-1}c)_{+} \\
&\ls 2C_2 k^p\sum^{\infty}_{n=m_k+1} (r^{n+1}+1)(2r-1)r^{-qn}
f(D_2(r,q)(r^q/2)^{(n-m_k-1)/p}) \\
&\ls 4C_2 r^2 k^p 
\sum_{n=m_k+1}^{\infty} r^{n(1-q)}f(D_2(r,q)(r^q/2)^{(n-m_k-1)/p}) \\
&\ls D_3(r,q)k^{\frac{p}{q}} c^{1-\frac{1}{q}} \sum^{\infty}_{n=m_k+1} r^{(n-m_k-1)(1-q)}
f(D_2(r,q)(r^q/2)^{-(n-m_k-1)/p}),
\end{align*}
where $D_3(r,q)=4C_2 r^{3-q} M_0^{\frac{1}{q}-1}$. Note that in the last line we have used that $c/M_0\ls k^pr^{-m_kq}$ which implies that
$$
k^p\ls M_0^{\frac{1}{q}-1} k^{\frac{p}{q}}c^{1-\frac{1}{q}}r^{-m_k(1-q)}.
$$
We have to consider
$$
D_4(r,q,p)=D_3(r,q)\sum^{\infty}_{l=0}r^{l(1-q)}f(D_2(r,q,p)(r^q/2)^{l/p})  
$$
but by our growth condition (\ref{nier:2}) $D_4(r,q,p)$ is finite and does not depend on $k$ nor $m_0$. Therefore 
we finally derive 
$$
\|TV^c(X)\|_p\ls (D_1(r,q)+D_2(r,q)D_4(r,q,p))k^{\frac{p}{q}}c^{1-\frac{1}{q}}.
$$
It ends the proof.
\end{dwd}
The consequence of the above theorem and our argument from
the beginning of Section \ref{sect2} is that if a process $X$ satisfies (\ref{nier:0}) for all $k\gs 1$
and $p\ls 1$ then it also satisfies (\ref{nier:1}) and hence by Theorem \ref{thm1}
its truncated variation satisfies the following concentration inequality. 
\begin{coro}\label{cor1}
Suppose that (\ref{nier:0}) is satisfied fpr all $k\gs 1$ and $p\ls 1$ then
$$
\PM(TV^c(X)\gs Duc^{1-\frac{1}{q}})\ls \exp(-u^{\frac{q}{p}}),\;\;\mbox{for}\;u\gs 1, 
$$
where $D$ is a universal constant.
\end{coro}
In particular Corollary \ref{cor1} works for any fractional Brownian motions $X_H$ with Hurst coefficient
$H\in (0,1)$. Indeed for each $k\gs 1$ the process $X_H$ satisfies
$$
\|X_H(t)-X_H(s)\|_k\ls Ck^{1/2}|t-s|^{H/2},
$$
so (\ref{nier:0}) holds for $p=\frac{H}{2}$. 
Therefore by Corollary \ref{cor1} we get
\begin{theo}\label{thm2}
For a fractional Brownian motion $X_H$ with Hurst coefficient $H\in (0,1)$ 
$$
\PM(TV^c(X_H)>Duc^{1-\frac{1}{H}})\ls e^{-u^{2H}},\;\;\mbox{for}\;u\gs 1,
$$
where $D$ is a universal constant.
\end{theo}


\begin{thebibliography}{}
%
%

\bibitem{Bed1} \textsc{Bednorz, W.} A theorem on majorizing measures, \textit{Annals of Probability}, 
{\bf 34}, no. 5, 1771-1781, (2006)

\bibitem{Bed2} \textsc{Bednorz, W.} On a Sobolev type inequality and its applications \textit{Studia Mathematica}, {\bf 176}, no. 2, 113-137, (2006) 

\bibitem{Bed3} \textsc{Bednorz, W.} (2010), Majorizing measures on metric spaces, \textit{C.R. math. Acad. Sci. Paris}, {\bf 348}, no. 1-2, 75-78, (2010)

\bibitem{Be-Lo} \textsc{Bednorz, W. and Lochowski, R.} (2015) Integrability and concentration of the truncated variation for the sample paths of fractional Brownian motions, diffusions and Lévy processes,
{\bf 21}, 437-464.

\bibitem{Li} \textsc{Li, W.V. and Shao Q.M.} (2001) Gaussian Processes, Small Ball Probabilities and Applications. \textit{Stoch. Proc.: Theory and methods   } Handbook of Statistics, Elsevier.

\bibitem{Fer1} \textsc{Fernique, X.}  Caract{\'e}risation de processus {\'a} trajectoires major{\'e}es ou continues.
(French) \textit{S{\'e}minaire de Probabilit{\'e}s XII.} Lecture Notes in Math. {\bf 649}, 691--706,
Springer, Berlin. (1978)

\bibitem{Fer2} \textsc{Fernique, X.}  R{\'e}gularit{\'e} de fonctions al{\'e}atoires non gaussiennes.
(French) \textit{{\'E}cole
d'{\'E}t{\'e} de Probabilit{\'e}s de Saint-Flour XI-1981.} Lecture Notes in Math. {\bf 976}, 1--74,
Springer, Berlin. (1983)

\bibitem{Kwap} \textsc{Kwapie\'n, S. and Rosi\'nski, J.} Sample H{\"o}lder continuity of stochastic processes and majorizing measures. \textit{Seminar on Stochastic Analysis, Random Fields and Applications IV, Progr. in Probab.} {\bf 58}, 155--163.
Birkh{\"a}user, Basel. (2004)

\bibitem{Le-Ta} \textsc{Ledoux, M. and Talagrand, M.}  Probability in Banach spaces.
Isoperimetry and processes. \textit{Results in Math. and Rel. Areas (3).} {\bf 23}, xii+480 pp.  Springer-Verlag, Berlin. (1991)

\bibitem{Loch}\textsc{\L ochowski, R.}  On the Generalisation of the Hahn-Jordan Decomposition for Real C\`{a}dl\`{a}g Functions. \textit{Colloq. Math.}, to appear.

\bibitem{Loch1}\textsc{\L ochowski, R.}  Truncated variation, upward truncated variation and downward truncated variation of Brownian motion with drift -- Their characteristics and applications. \textit{Stochastic Proc. Appl.} {\bf 121}, no. 2, 378-393 (2011)

\bibitem{Loch2}\textsc{\L ochowski, R.} Pathwise stochastic integration with finite variation processes uniformly approximating c\`{a}dl\`{a}g processes. \textit{arXiv} (2013)

\bibitem{LM}\textsc{\L ochowski, R., Mi\l o\'{s}, P.}  On truncated variation, upward truncated variation and downward truncated variation for diffusions. \textit{Stochastic Proc. Appl.} {\bf 123}, no. 2, 446-474 (2013)

\bibitem{P08} Picard, Jean A tree approach to p-variation and to
integration. Ann. Probab. 36, no. 6, 2235\textendash{}2279 (2008) 

\bibitem{RY} \textsc{Revuz, D. and Yor, M.} Continuous martingales and Brownian motion,
volume 293 of Grundlehren der Mathematischen Wissenschaften [Fundamental Principles of Mathematical Sciences], Springer, (2005)

\bibitem{Sato:1999}
\textsc{ Sato, K.-I.}, L\'{e}vy Processes and Infinite Divisibility. 
\textit{Cambridge Univ. Press}. (1999)

\bibitem{Tal1} \textsc{Talagrand, M.} Sample boundedness of stochastic processes under increment conditions. \textit{Annals of Probability} {\bf 18}, no. 1, 1-49 (1990)

\bibitem{Tal2} \textsc{Talagrand, M.} The generic chaining. \textit{Springer-Verlag}. (2005)



\end{thebibliography}
\end{document}